\documentclass{article}

\usepackage{hyperref}
\usepackage[bibtex-style]{amsrefs}
\usepackage{pstricks}
\usepackage{amsmath, amsthm, amssymb}

\newtheorem{theorem}{Theorem}
\newtheorem{corollary}[theorem]{Corollary}

\newtheorem{proposition}[theorem]{Proposition}

\newcommand{\C}{{\mathcal C}}
\renewcommand{\S}{{\mathcal S}}

\newcommand{\involvedIn}{\preceq}
\newcommand{\doublyDom}{\sqsubseteq}

\newcommand{\dominanceOrder}{\trianglelefteq}
\newcommand{\shape}{\mathop{\mathrm{sh}}}
\newcommand{\avoids}{\mathop{\mathrm{Av}}}

\title{Young Classes of Permutations}
\author{Michael Albert \\%
Department of Computer Science, University of Otago, Dunedin,
New Zealand \\
\url{malbert@cs.otago.ac.nz}
}

\begin{document}
\maketitle
\begin{abstract}
  We characterise those classes of permutations closed under taking
  subsequences and reindexing that also  have the property
  that for every tableau shape either every permutation of that shape
  or no permutation of that shape belongs to the class. The
  characterisation is in terms of the dominance order for partitions
  (and their conjugates) and shows that for any such class there is a
  constant $k$ such that no permutation in the class can contain both
  an increasing and a decreasing sequence of length $k$.
\end{abstract}

\section{Introduction}

The most well known notion of the \emph{shape} of a permutation arises
from the Robinson-Schensted-Knuth correspondence, whereby the shape of
$\pi \in \S_n$ is taken to be the shape of either of the Young
tableaux produced from it by the RSK-algorithm. On the other hand, we
might think of the shape of $\pi$ as being described by its graph, i.e.~the set of points $\{(i, \pi(i)) \, | \, 1 \leq i \leq n \}$. 

In the first context, subshapes can be defined in terms of various
natural orderings on the set of partitions. In the second context, there is a natural notion of \emph{subpermutation} which we can introduce. The formal definition is spelled out in the next section but informally, $\sigma \in \S_k$ is a subpermutation of $\pi \in \S_n$ if, there is some $k$ element subset of $\{(i, \pi(i)) \, | \, 1 \leq i \leq n \}$ whose relative ordering with respect to both axes is the same as that of $\{(j, \sigma(j)) \, | \, 1 \leq j \leq k \}$. 

It seems natural to investigate the manner in which these two concepts
interact. Unfortunately, the basic answer is: not very well. In
particular Adin and Roichman \cite{AR02} began such an investigation
and reported indicative counterexamples to some tempting false
conjectures, as well as more positive results relating to rectangular
and hook shapes. More recently, at the Permutation Patterns 2010
conference at Dartmouth College, Panova \cite{Panova10}, \cite{Crites10} and
Tiefenbruck \cite{Tiefenbruck10} reported on results in this area when
attention is restricted to permutations avoiding certain specific
patterns.

Since the two concepts seem somewhat orthogonal with respect to \emph{specific} permutations, it remains to consider whether or not they interact well with respect to \emph{collections} of permutations. Indeed Schensted's theorem (and Greene's generalisation of it) establish a family of such relationships. For example: the collection of all permutations not containing a decreasing sequence of length 3 is the same as the collection of all permutations whose tableaux have at most two rows.

A set of permutations closed under taking subpermutations is
called a \emph{permutation class} (or simply \emph{class}). Motivated by the result mentioned above arising from Schensted's theorem, 
we investigate the question:
\begin{quote}
Which classes, $\C$, of permutations have the property that if $\pi$ and
$\pi'$ are two permutations of the same shape then either both, or
neither belong to $\C$?
\end{quote}
Such classes are precisely those closed in the natural equivalence
relation on permutations of ``having the same shape''. Classes of
this type will be called \emph{Young classes}. 

Greene's
generalisation of Schensted's theorem provides a natural candidate for
a condition on the set of shapes that may occur in a Young class, and
the main result of this paper is to show that this condition is both
necessary (which is the new result) and sufficient (which is essentially a corollary of the generalisation). The most notable
consequence of our characterisation of these classes is that for any
proper Young class there is a constant $k$ such that no permutation in
the class contains both increasing and decreasing subsequences of
length $k$.

\section{Definitions}

We refer the reader to either \cite{Sagan01} or \cite{Stanley99} for
more detailed discussion of the RSK-correspondence and related
matters. In this section we will collect only the basic definitions
and results needed for the remainder of this paper.

A \emph{partition}, $\lambda$, is a weakly decreasing sequence of
positive integers, i.e. $\lambda = (\lambda_1, \lambda_2, \ldots,
\lambda_k)$ where $\lambda_1 \geq \lambda_2 \geq \cdots \geq \lambda_k
> 0$. We say that $\lambda$ is a partition of $n$, and write
$|\lambda| = n$ where $n = \sum_i
\lambda_i$. The individual values $\lambda_i$ are called the
\emph{parts} of $\lambda$. When writing a partition as a sequence, repeated values
may be represented by exponentiation. A partition is often represented by its \emph{Young
  diagram}, an array of boxes of the appropriate lengths (we use the
English or ``top down'' convention). For instance:

\centerline{
\psset{xunit=0.3}
\psset{yunit=0.3}
\begin{pspicture}(-11,-1)(4,4)
\rput(-7,2){$(3,3,2,1) = (3^2, 2, 1) = $}
\psline(0,0)(1,0)
\psline(0,1)(2,1)
\psline(0,2)(3,2)
\psline(0,3)(3,3)
\psline(0,4)(3,4)
\psline(0,0)(0,4)
\psline(1,0)(1,4)
\psline(2,1)(2,4)
\psline(3,2)(3,4)
\end{pspicture}
}

It is sometimes convenient to pretend that a partition may be
extended by (any number of) trailing $0$'s. The
\emph{conjugate partition}, $\lambda^*$ of a partition $\lambda$ is the sequence
given by:
\[
\lambda^*_i = \left| \{ j \, : \, \lambda_j \geq i \} \right|.
\]
The \emph{sum} $\lambda + \mu$ of two partitions is defined in the
obvious way (its $i^{\mbox{\scriptsize th}}$ element is $\lambda_i +
\mu_i$), and their \emph{conjugate sum} is defined by: 
\[ 
\lambda +^{*} \mu = \left( \lambda^* + \mu^* \right)^*.
\]
On Young diagrams these operations correspond to concatenating the two
partitions horizontally or vertically, and then shoving blocks
leftwards or upwards to eliminate gaps.

If $\lambda$ and $\mu$ are partitions, we say that $\mu$
\emph{dominates} $\lambda$ and write $\lambda \dominanceOrder \mu$ if,
for all $k \geq 1$:
\[
\sum_{i=1}^k \lambda_i \leq \sum_{i=1}^k \mu_i.
\]
If both $\lambda \dominanceOrder \mu$ and $\lambda^* \dominanceOrder
\mu^*$ then we say that $\mu$ \emph{doubly dominates} $\lambda$ and
write $\lambda \doublyDom \mu$. Since conjugation reverses the dominance order on partitions of $n$ , if $|\lambda| =
|\mu|$ and $\lambda \doublyDom \mu$ then $\lambda = \mu$.

The \emph{shape} of permutation $\pi \in \S_n$ is the partition
$\lambda$ of $n$ obtained when constructing its Young tableaux via the
RSK-correspondence. In this case we write $\shape(\pi) = \lambda$. We also say that $n$ is the \emph{length} of $\pi$ since we frequently consider $\pi$ and the sequence of its values interchangeably.

A permutation $\sigma$ \emph{is involved in} $\pi$ if the terms of
some subsequence of the values of $\pi$ occur in the same relative order as the sequence of values
of $\sigma$. In this case we write $\sigma \involvedIn \pi$, and also
say that $\sigma$ is a \emph{subpermutation} of $\pi$. For example,
$231 \involvedIn 421635$ because of the subsequence $463$. A
\emph{permutation class} is a set of permutations closed downwards
under $\involvedIn$. A permutation class can also be described as the
set of permutations that do not involve (or \emph{avoid}) some given
set of permutations $X$. In that case we write $\C = \avoids(X)$. 

A \emph{Young class} is a class of permutations such that for each
partition $\lambda$ either all the permutations of shape $\lambda$ or
none of them belong to the class. The prototypical examples of Young
classes are the collection of all permutations not including an
increasing subsequence of some fixed length $a$. By Schensted's
theorem \cite{Schensted61}, the shapes of all
such permutations have every part smaller than $a$, and the converse
also holds.

The use of Greek letters for both partitions and permutations is
unfortunately standard. We will attempt to minimise the resulting
confusion by reserving $\lambda$ and $\mu$ for partitions.

\section{Greene's Theorem and double domination}

Greene's generalisation \cite{Greene74} (or see Theorem 3.5.3 of
\cite{Sagan01}) of Schensted's theorem
provides a link between the shape of a permutation
$\pi$ and the involvement of certain patterns in $\pi$.

\begin{theorem}
For any permutation $\pi$ of tableau shape $\lambda$, and any $k$, the
longest subpermutation of $\pi$ that can be written as a union of
$k$ increasing subsequences has length $\lambda_1 + \lambda_2 + \cdots
+ \lambda_k$, and the longest subpermutation of $\pi$ that can be
written as a union of $k$ decreasing subsequences has length
$\lambda^\ast_1 + \lambda^\ast_2 + \cdots + \lambda^\ast_k$.
\end{theorem}

The increasing subsequences that make up the subpermutations whose existence is guaranteed by Greene's theorem can of course be taken to be disjoint, but they cannot in general be taken to have lengths $\lambda_1$, $\lambda_2$, \dots, $\lambda_k$. A simple counterexample is the permutation $348951267$ where $\lambda_1 = 5$, and $\lambda_2 = 3$. However, the unique increasing subsequence of length 5 is 34567, and the remaining elements, 8912, contain no three element increasing sequence. On the other hand, we can find two disjoint four element increasing sequences: 3489 and 1267.

Note that the permutations that avoid
$(k+1) k \cdots 321$ are precisely those which can be written as a union of
$k$ increasing subsequences. So $\lambda_1 + \lambda_2 +
\cdots + \lambda_k$ is the maximum length of an element of
$\avoids((k+1) k \cdots 321)$ which is involved in $\pi$. Of course a
dual statement applies with reference to the second part of the
theorem and the class $\avoids(123\cdots (k+1))$.

The following corollary to Greene's theorem is immediate:

\begin{corollary}
\label{COR_ShapeInvolvement}
Let $\sigma$ and $\pi$ be permutations. If $\sigma \involvedIn \pi$
then $\shape(\sigma) \doublyDom \shape(\pi)$.
\end{corollary}
\begin{proof}
If $\sigma \involvedIn \pi$ then, for any $k$
the longest subsequence of $\pi$ which can be written as the union of
$k$ increasing (decreasing) subsequences is at least as long as the
corresponding subsequence of $\sigma$, since that subsequence occurs
in $\pi$.
\end{proof}

A strong converse to the corollary is easily seen to be false. For example
$(2,2) \doublyDom (3,1,1)$ but the permutation $54123$ of shape
$(3,1,1)$ has no subpermutation of shape $(2,2)$. However, even the
weaker converse ``if $\lambda \doublyDom \mu$ then for some
permutations $\sigma$ and $\pi$ with $\shape(\sigma) = \lambda$ and
$\shape(\pi) = \mu$, $\sigma \involvedIn \pi$'' is also false. For, take
$\lambda = (2,2,2)$ and $\mu = (4,1,1,1,1)$. If $\sigma$ has shape
$\lambda$, then to extend it to a permutation having a longest
increasing subsequence of length $4$ a pair of elements in increasing
order must be added. However, to extend it to a permutation whose
longest decreasing subsequence has length $5$ a pair of elements in
decreasing order must be added. To accomplish both therefore requires
adding at least three elements, and thus rules out having shape $\mu$.
On the other hand, we will see below that if $\lambda \doublyDom \mu$
and $| \mu | = | \lambda | + 1$ then there are such permutations
$\sigma$ and $\pi$. This will be the central result that leads to a
characterisation of Young classes.

For our subsequent arguments we will need to characterize the covering pairs for the relation $\doublyDom$ on partitions. The main point is to show that only the obvious pairs (those whose sizes differ by 1) are covering pairs.

\begin{proposition}
\label{PROP_InvolvementCovers}
Suppose that $\lambda \doublyDom \mu$, and $|\mu| - |\lambda| >
1$. Then there exists $\mu'$ with $\lambda \doublyDom \mu'
\doublyDom \mu$ and $|\mu| - |\mu'| = 1$.
\end{proposition}

\begin{proof}
  We begin by eliminating some trivial or irrelevant cases.
  \begin{itemize}
  \item
  If $\lambda_1 = \mu_1$ or $\lambda^*_1 = \mu^*_1$, then the result follows inductively by ignoring the first parts of $\lambda$ and $\mu$ (or of their conjugates).
  \item
  If, for all $r$, $\sum_{i=1}^r \mu_i > \sum_{i=1}^r \lambda_i$, then we can form $\mu'$ by deleting the cell of the lowest row of the rightmost column of $\mu$. This decreases each sum of row lengths by at most 1 (compared to $\mu$) and leaves all but the final sum of column lengths unchanged. Since $\mu_1 > \lambda_1$ the only relevant comparisons of column sums occur before the final column, and so $\lambda \doublyDom \mu'$.
  \item
  A similar argument applies if, for all $c$, $\sum_{i=1}^c \mu^*_i > \sum_{i=1}^c \lambda^*_i$.
\end{itemize}

Assuming that none of the cases above hold, we are left with the following conditions: $\mu_1 > \lambda_1$; $\mu^*_1 > \lambda^*_1$; for some $r > 1$, $\sum_{i=1}^r  \mu_i = \sum_{i=1}^r \lambda_i$; and for some $c > 1$, $\sum_{i=1}^c \mu^*_i = \sum_{i=1}^c \lambda^*_i$. Choose the least possible values of $r$ and $c$ such that these equalities hold, and observe that $\mu_r < \lambda_r$ and $\mu^*_c < \lambda^*_c$.

  Suppose that $\mu_r \geq c$ and $\mu^*_c \geq r$, i.e. that $\mu$
  contains the upper left $r \times c$ rectangle. Since $\lambda_r >
  \mu_r$ and $\lambda^*_c > \mu^*_c$, $\lambda$ also contains the upper
  left $r \times c$ rectangle. Divide the partitions $\lambda$ and $\mu$ up into four regions as shown in Figure~\ref{FIG_Subdiv}.

\begin{figure}
\centerline{
\psset{xunit=0.2}
\psset{yunit=0.2}
\begin{pspicture}(-2,-2)(25,25)
   \psline(0,0)(24,0)(24,24)(0,24)(0,0)
   \psline(12,0)(12,24)
   \psline(0,16)(24,16)
\rput(6,20){$A$}
\rput(18,20){$B$}
\rput(6,10){$C$}
\rput(18,10){$D$}
\rput(-1,16){$r$}
\rput(12,-1){$c$}
\end{pspicture}
}
\caption{The subdivision of $\pi$ and $\mu$ used in the proof of Proposition \ref{PROP_InvolvementCovers}.}
\label{FIG_Subdiv}
\end{figure}
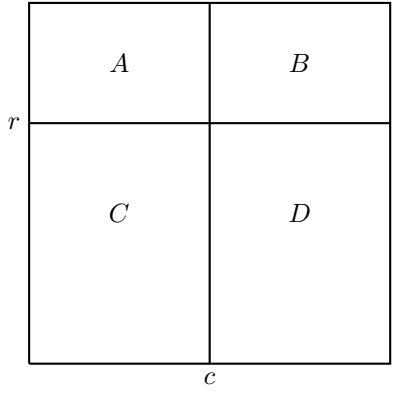

  Both partitions contain the same number of cells in $A \cup B$ (since their $r^{\mbox{\scriptsize th}}$ row sums agree) and also in $A \cup C$. Since both contain the entire region $A$, they contain the same number of cells in $B$ and in $C$. So, there are at least two more cells of $\mu$ in region $D$ than of $\lambda$, and the portion of $\mu$ in $D$ doubly dominates the corresponding portion of $\lambda$. The result for this case now follows inductively by using it in $D$.

  So we may now additionally assume that either $\mu_r < c$ or $\mu^*_c < r$ and
  hence that there is at least one square in the upper left $r \times
  c$ rectangle not belonging to $\mu$. Form $\mu'$ by removing a cell
  from the rightmost column, and a cell from the bottommost row of
  $\mu$ and adding a cell within the upper left $r \times c$
  rectangle. Then certainly $\mu' \doublyDom \mu$ and we claim that
  $\lambda \doublyDom \mu'$. Consider first the row sums of
  $\mu'$. Some among the first $r-1$ may be one smaller than the
  corresponding sums of $\mu$ but, by the choice of $r$ they are still
  at least as large as those of $\lambda$. From rows $r$ through to
  the final row they equal those of $\mu$, and in the final row are
  again one smaller. However as $|\mu| > |\lambda|$ they are all at
  least as great as those of $\lambda$. The argument that the column
  sums of $\mu'$ dominate those of $\lambda$ is exactly similar.
\end{proof}

We now know that if $\lambda \doublyDom \mu$ is a covering pair then $|\mu| = |\lambda| + 1$. We can also bound the difference between row sums of $\lambda$ and of $\mu$. Specifically:

\begin{proposition}
Suppose that $\lambda \doublyDom \mu$ and $|\mu| = |\lambda|+1$. Then, for all $r$, 
\[
0 \leq \sum_{i=1}^r \mu_i - \sum_{i=1}^r \lambda_i \leq 1.
\]
\end{proposition}

\begin{proof}
The left hand inequality is part of the definition of $\lambda \doublyDom \mu$. For the remainder, suppose that it failed. Choose the least $r$ such that $\sum_{i=1}^r \mu_i - \sum_{i=1}^r \lambda_i \geq 2$. Then $\mu_r > \lambda_r$. Let $c = \lambda_r$. Consider again the diagram in the preceding proof. Region $A$ is fully occupied in both $\lambda$ and $\mu$, while region $B$ contains at least two more cells of $\mu$ than of $\lambda$, and region $D$ contains no cells of $\lambda$. So region $C$ contains at least one more cell of $\lambda$ than of $\mu$ but this (along with $A$ being fully occupied) implies $\sum_{i=1}^c \lambda_i^* > \sum_{i=1}^c \mu_i^* $ and so contradicts $\lambda \doublyDom \mu$.
\end{proof}
 
\section{Young classes}

To characterise Young classes we need to introduce two more operations
on permutations, and show how they affect the shape of the resulting
tableaux. Let permutations $\alpha$ and $\beta$ of lengths $n$ and $k$
be given. Define:
\begin{eqnarray*}
\alpha \oplus \beta &=& (\alpha_1, \ldots, \alpha_n, n + \beta_1,
\ldots, n + \beta_k) \\
\alpha \ominus \beta &=& (k+\alpha_1, \ldots, k+\alpha_n, \beta_1,
\ldots, \beta_k).
\end{eqnarray*}
Informally, $\alpha \oplus \beta$ stacks $\beta$ above and to the
right of $\alpha$, while $\alpha \ominus \beta$ stacks $\beta$ below and
to the right of $\alpha$.

\begin{proposition}
The operations $\oplus$ and $\ominus$ affect shape as follows:
\begin{eqnarray*}
\shape(\alpha \oplus \beta) &=& \shape(\alpha) + \shape(\beta) \\
\shape(\alpha \ominus \beta) &=& \shape(\alpha) +^* \shape(\beta).
\end{eqnarray*}
\end{proposition}

\begin{proof}

  The proof of both parts is most easily seen by considering Greene's
  theorem. In $\alpha \oplus \beta$ a subsequence that is the union of
  $k$ increasing sequences can be divided into such a subsequence of
  $\alpha$ and such a subsequence of $\beta$ (and the union of any two
  such subsequences is again of the same type). So, the sum of the
  first $k$ parts of $\shape(\alpha \oplus \beta)$ is equal to the sum
  of the first $k$ parts of $\shape(\alpha)$ plus the sum of the first
  $k$ parts of $\shape(\beta)$. The result follows immediately. The
  corresponding result for $\ominus$ likewise follows by considering
  unions of decreasing subsequences.

\end{proof}

We remark that the result for $\oplus$ is also easily seen by
considering the actual bumping operations of the RSK algorithm, but
the result for $\ominus$ is not so transparent in that context (though it is true that the 
operations of RSK ``push straight down columns'').

\begin{proposition}
\label{PROP_DDCoverIsInvolvementWitnessed}
Let $\lambda \doublyDom \mu$ be a covering pair for the double domination ordering. Then there exist permutations $\sigma$ and $\pi$ such that: $\shape(\sigma) = \lambda$, $\shape(\pi) = \mu$ and $\sigma \involvedIn \pi$.
\end{proposition}

\begin{proof}
The proof is by induction on $|\lambda|$. The base case $|\lambda| = 0$ is completely trivial. Consider now a covering pair $\lambda \doublyDom \mu$. 

Suppose first that for some $i$, $\lambda_i = \mu_i = k$, and choose $i$ to be the least such. Let $\lambda'$ and $\mu'$ be obtained by deleting the $i^{\mbox{\scriptsize th}}$ part from $\lambda$ and $\mu$ respectively. Inductively, choose $\sigma' \involvedIn \pi'$ with $\shape(\sigma') = \lambda'$ and $\shape(\pi') = \shape(\sigma')$. Let $\iota = 123 \cdots k$, and let $\sigma = \iota \ominus \sigma'$ and $\pi = \iota \ominus \pi'$. Then $\sigma \involvedIn \pi$, and 
\[
\begin{array}{c}
\shape(\sigma) = (k) +^*  \lambda' = \lambda \\
\shape(\pi) = (k) +^*  \mu' = \mu .
\end{array}
\]
Thus the result is true in this case, and by a parallel argument in the case where $\lambda_i^* = \mu_i^*$ for some $i$.

So we now assume that no part of $\lambda$ is the same size as the corresponding part of $\mu$ (and that this also holds for their conjugates). Since $\lambda \doublyDom \mu$ and since we know that the row sums of $\mu$ can exceed the corresponding row sums of $\lambda$ by at most 1, the row sums of $\mu$ must alternately exceed by 1 and equal the corresponding row sums of $\lambda$ (otherwise we would have two equal parts). So:
\begin{eqnarray*}
\mu_1 &=& \lambda_1 + 1 \\
\mu_2 &=& \lambda_2 - 1 \\
\mu_3 &=& \lambda_3 + 1 \\
{} &\cdots& 
\end{eqnarray*}
The partition $\mu$ must have exactly one more part than $\lambda$ and this final part must be 1. All together this implies that for some positive integer $k$: 
\begin{eqnarray*}
\mu &=& (2k+1, 2k-1, 2k-1, 2k-3, 2k-3, \dots, 3, 3, 1) \\
\lambda &=& (2k, 2k, 2k-2, 2k-2, \dots, 2, 2).
\end{eqnarray*}
Since we shall need to do an induction on $k$, denote the first of these by $\mu^{(k)}$ and the second by $\lambda^{(k)}$. It remains to prove the result in precisely this case. The general construction is based on the first case of this sequence $(2, 2) \doublyDom (3,1)$. Here we may take $\sigma = 2413$ and $\pi = 25314$. For use in the general case, define:
\[
\theta^{(k)} = 1 \oplus 21 \oplus 321 \oplus \dots \oplus k (k-1) \dots 321.
\]
The key property of $\theta^{(k)}$ is that its maximal increasing and decreasing subsequences both have length $k$ and the deletion of either a maximal increasing or a maximal decreasing subsequence leaves a copy of $\theta^{(k-1)}$.

Let $\sigma^{(k)}$ be formed from $\sigma = 2413$ by ``inflating'' each of the elements into a copy of $\theta^k$, and correspondingly let $\pi^{(k)}$ be obtained from $\pi = 25314$ by inflating the elements 2, 5, 1, and 4 in the same way, while leaving the point corresponding to 3 as a singleton. For instance, $\pi^{(3)} \in \S_{25}$ is shown in Figure~\ref{FIG_Pi3} (and $\sigma^{(3)}$ is obtained from it by deleting the central point).

\begin{figure}
\newcommand{\thet}{%
\gsave
\psline(0,0)(7,0)(7,7)(0,7)(0,0)
\stroke
\grestore
\gsave
\pscircle(1,1){2pt}
\fill
\grestore
\gsave
\pscircle(2,3){2pt}
\fill
\grestore
\gsave
\pscircle(3,2){2pt}
\fill
\grestore
\gsave
\pscircle(4,6){2pt}
\fill
\grestore
\gsave
\pscircle(5,5){2pt}
\fill
\grestore
\gsave
\pscircle(6,4){2pt}
\fill
\grestore
}

\centerline{
\psset{xunit=0.2}
\psset{yunit=0.2}
\begin{pspicture}(-2,-2)(33,33)
   \pscustom[fillstyle=solid, fillcolor=black]{%
   \translate(0,7)
   \thet
   }
   \pscustom[fillstyle=solid, fillcolor=black]{%
   \translate(7,24)
   \thet
   }
\pscustom[fillstyle=solid, fillcolor=black]{%
   \translate(17,0)
   \thet
   }
\pscustom[fillstyle=solid, fillcolor=black]{%
   \translate(24,17)
   \thet
   }
   \pscircle*[fillstyle=solid, fillcolor=black](15.5,15.5){3pt}
\end{pspicture}
}
\caption{The permutation $\pi^{(3)}$ illustrating the proof of Proposition \ref{PROP_DDCoverIsInvolvementWitnessed}.}
\label{FIG_Pi3}
\end{figure}
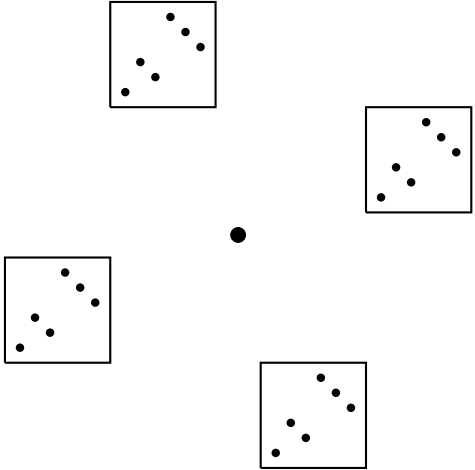
We claim that $\lambda^{(k)} = \shape(\sigma^{(k)})$ and $\mu^{(k)} = \shape(\pi^{(k)})$ (and this of course will finish the proof). This is proved by induction on $k$, with the case $k = 1$ already having been considered above. 

We will provide the argument for $\pi^{(k)}$ in detail. Clearly an increasing sequence in $\pi^{(k)}$ can contain only elements in two of the inflated regions, together with possibly the central point (if the regions used are the leftmost and rightmost). In particular, the maximum increasing sequence has length $2k+1$. However, if two increasing sequences include one using the central point then they can use only three of the four inflated regions, giving a maximum possible length of $4k$ (in two of the three regions we can only use an increasing sequence so obtaining $2k$ elements, in one we might use two increasing sequences giving a further $2k-1$ elements, and finally we could include the central point). On the other hand we could also achieve a length of $4k$ by choosing an increasing sequence of length $2k$  from the two leftmost regions, and another one from the two rightmost regions. So, the first two rows of $\shape(\pi^{(k)})$ have lengths $2k+1$ and $2k-1$. In our second construction we found two increasing sequences in $\pi^{(k)}$ of total length $4k$ whose removal leaves a copy of $\pi^{(k-1)}$. By induction, the row sums of $\shape(\pi^{(k)})$ are at least as great as those of $\mu^{(k)}$. The same argument applied to decreasing sequences gives the corresponding result for column sums. So $\mu^{(k)} \doublyDom \shape(\pi^{(k)})$. Therefore, as the two partitions have the same size, they must be equal.

The argument for $\sigma^{(k)}$ is entirely similar, and in fact simplified by the omission of the central point.
\end{proof}

Finally, we can prove our main result:

\begin{theorem}
\label{THM_YoungClasses}
A permutation class $\C$ is a Young class if and only if there is some
downwards closed set $D$ in the double domination order of partitions and
$\C = \{ \pi \, : \, \shape(\pi) \in D\}$.
\end{theorem}

\begin{proof}
  By definition, $\C$ must be a union of shape equivalence classes,
  over a set $D$ of shapes. To see that it is necessary that $D$ be a
  downwards closed set in the double domination order suppose that $\mu \in D$. By
  Proposition~\ref{PROP_DDCoverIsInvolvementWitnessed}, if $\lambda \doublyDom \mu$ has
  $|\lambda| + 1 = |\mu|$, there are two permutations $\sigma$
  and $\pi$ with $\shape(\sigma) = \lambda$, $\shape(\pi) = \mu$ and
  $\sigma \involvedIn \pi$. As $\pi \in \C$ we must have $\sigma \in
  \C$ and hence $\lambda \in D$. However, by
  Proposition~\ref{PROP_InvolvementCovers}, the
  full double domination relation is the transitive closure of these
  single element covers. So, $D$ must be downwards closed.

  On the other hand, Corollary~\ref{COR_ShapeInvolvement} says
  precisely that if $D$ is downwards closed, then so is $\C$.
\end{proof}

Note that for any partition $\lambda$ there exists a hook, i.e.~a
partition of the form $(n, 1^m)$, such that $\lambda \doublyDom (n,
1^m)$ (essentially we can replace each element not in the first row or
column by a pair of elements, one in the first row and one in the
first column). It follows that in any proper Young class there are
bounds, say $a$ and $d$ such that every permutation in the class
either has no increasing subsequence of length $a$ or no decreasing
subsequence of length $d$. So every proper Young class is contained in
the union of two of the prototypical Young classes already described
by Schensted's theorem. Since, at the cost of enlarging one of those two classes, we could always take the monotone permutations which they avoid to have equal length, the final remark of the introduction is justified.

\section{Conclusions and Open Problems}

We have considered classes of permutations closed under the equivalence relation of having tableaux of the same shape. The general
relationship between the involvement relationship for permutations and
their shapes seems to be somewhat obscure. However, by restricting
attention to these Young classes we were able to demonstrate a close
connection with the double domination order for partitions.

More detailed connections (e.g. to describe the bases of these
classes) seem to be difficult to achieve except in cases closely
related to the prototypical Young classes provided by Schensted's theorem.

Returning to the vexing case of the relationship between involvement and shape for pairs of permutations (or tableau shapes), the most obvious remaining open question is:
\begin{quote}
Give necessary and sufficient conditions on partitions $\lambda
\doublyDom \mu$ which guarantee the existence of permutations
$\sigma$ and $\pi$ with $\shape(\sigma) = \lambda$, $\shape(\pi) =
\mu$ and $\sigma \involvedIn \pi$.
\end{quote}
The example following Corollary~\ref{COR_ShapeInvolvement} would seem
to suggest that such conditions would involve a strengthening or
refinement of Greene's theorem, but the details remain elusive.

\section{Acknowledgments}

The author would like to thank the University of Otago Theory Group
(Mike Atkinson, Robert Aldred, Dennis McCaughan) who inspired the
preliminary stages of this work in early 2009 and the organisers of the
Permutation Patterns 2010 conference.

The author would also like to thank two careful and thorough referees who exposed an oversight in the main argument in an earlier version of this paper.


\begin{bibdiv}
\begin{biblist}

\bib{AR02}{article}{
   author={Adin, Ron M.},
   author={Roichman, Yuval},
   title={Shape avoiding permutations},
   journal={J. Combin. Theory Ser. A},
   volume={97},
   date={2002},
   number={1},
   pages={162--176},
   issn={0097-3165},
   review={\MR{1879132 (2002j:05001)}},
   doi={10.1006/jcta.2001.3202},
}

\bib{Crites10}{article}{

  author={Crites, Andrew},
  author={Panova, Greta},
  author={Warrington, Gregory~S.},
  title={Separable permutations and {G}reene's theorem},
  date={2010},
  note={\url{http://arxiv.org/abs/1011.5491}}

}

\bib{Greene74}{article}{
   author={Greene, Curtis},
   title={An extension of Schensted's theorem},
   journal={Advances in Math.},
   volume={14},
   date={1974},
   pages={254--265},
   issn={0001-8708},
   review={\MR{0354395 (50 \#6874)}},
}

\bib{Panova10}{article}{
  author={Panova, Greta},
  title = {Separable permutations, Robinson-Schensted and shortest
    containing supersequences},
    conference={
    title={Permutation Patterns},
    address={Dartmouth College},
    date={2010}
    },
  eprint={http://www.math.dartmouth.edu/~pp2010/abstracts/Panova.pdf},
}


\bib{Sagan01}{book}{
   author={Sagan, Bruce E.},
   title={The symmetric group},
   series={Graduate Texts in Mathematics},
   volume={203},
   edition={2},
    publisher={Springer-Verlag},
   place={New York},
   date={2001},
   pages={xvi+238},
   isbn={0-387-95067-2},
   review={\MR{1824028 (2001m:05261)}},
}

\bib{Schensted61}{article}{
   author={Schensted, C.},
   title={Longest increasing and decreasing subsequences},
   journal={Canad. J. Math.},
   volume={13},
   date={1961},
   pages={179--191},
   issn={0008-414X},
   review={\MR{0121305 (22 \#12047)}},
}

\bib{Stanley99}{book}{
   author={Stanley, Richard P.},
   title={Enumerative combinatorics. Vol. 2},
   series={Cambridge Studies in Advanced Mathematics},
   volume={62},
    publisher={Cambridge University Press},
   place={Cambridge},
   date={1999},
   pages={xii+581},
   isbn={0-521-56069-1},
   isbn={0-521-78987-7},
   review={\MR{1676282 (2000k:05026)}},
}
		
\bib{Tiefenbruck10}{article}{
  author={Tiefenbruck, Mark},
  title = {231-avoiding permutations and the Schensted correspondence},
  conference={
    title={Permutation Patterns},
    address={Dartmouth College},
    date={2010}
    },
  eprint={http://www.math.dartmouth.edu/~pp2010/abstracts/Tiefenbruck.pdf},
}
	
		

	
\end{biblist}
\end{bibdiv}
\end{document}